\documentstyle[12pt,amssymb]{amsart} 
\def\SBIMSMark#1#2#3{
 \font\SBF=cmss10 at 10 true pt
 \font\SBI=cmssi10 at 10 true pt
 \setbox0=\hbox{\SBF Stony Brook IMS Preprint \##1}
 \setbox2=\hbox to \wd0{\hfil \SBI #2}
 \setbox4=\hbox to \wd0{\hfil \SBI #3}
 \setbox6=\hbox to \wd0{\hss
             \vbox{\hsize=\wd0 \parskip=0pt \baselineskip=10 true pt
                   \copy0 \break%
                   \copy2 \break%
                   \copy4 \break}}
 \dimen0=\ht6   \advance\dimen0 by \vsize \advance\dimen0 by 8 true pt
                \advance\dimen0 by -\pagetotal
 \dimen2=\hsize \advance\dimen2 by .25 true in
  \openin2=publishd.tex
  \ifeof2\setbox0=\hbox to 0pt{}
  \else 
     \setbox0=\hbox to 3.1 true in{
                \vbox to \ht6{\hsize=3 true in \parskip=0pt  \noindent  
                \input publishd.tex 
                \vfill}}
  \fi
  \closein2
  \ht0=0pt \dp0=0pt
 \ht6=0pt \dp6=0pt
 \setbox8=\vbox to \dimen0{\vfill \hbox to \dimen2{\copy0 \hss \copy6}}
 \ht8=0pt \dp8=0pt \wd8=0pt
 \copy8
 \message{*** Stony Brook IMS Preprint #1, #2 ***}
}

\newcommand{\iso}{\stackrel{\simeq}{\longrightarrow}}

\newcommand{\tree}{[x,y,\ldots,z]} 

\renewcommand{\marginpar}[1]{}
\catcode`\@=12

\def\Empty{}
\newcommand\oplabel[1]{
  \def\OpArg{#1} \ifx \OpArg\Empty {} \else
  	\label{#1}
  \fi}
		
\long\def\realfig#1#2#3#4{
\begin{figure}[htp]
\centerline{\psfig{figure=#2,width=#4}}
\caption[#1]{#3}
\oplabel{#1}
\end{figure}}

\newcommand{\comm}[1]{}

\input{psfig} 
\newtheorem{thm}{Theorem}
\newtheorem{cor}{Corollary}
\newtheorem{lem}{Lemma}

\newcommand{\thmref}[1]{Theorem~\ref{#1}}

\newcommand{\lemref}[1]{Lemma~\ref{#1}}
\newcommand{\corref}[1]{Corollary~\ref{#1}} 
\newcommand{\figref}[1]{Fig.~\ref{#1}}

\newcommand{\BBB}[1]{{\Bbb #1}}

\begin{document}
\SBIMSMark{1998/1a}{January 1998}{}

\title{Biaccessibility in quadratic Julia sets I: The locally-connected case} 
\author{Saeed Zakeri} 
\address{Department of Mathematics, SUNY at Stony Brook, NY 11794}
\email{zakeri@math.sunysb.edu}

\pagestyle{myheadings}
\markboth{\sc S. Zakeri}{\sc Biaccessibily in Quadratic Julia Sets I}
\begin{abstract}
Let $f:z\mapsto z^2+c$ be a quadratic polynomial whose Julia set $J$ is
locally-connected. We prove that the Brolin measure of the set of
biaccessible points in $J$ is zero except when $f(z)=z^2-2$ is the Chebyshev
quadratic polynomial for which the corresponding measure is one.  
\end{abstract} 
\maketitle

\noindent 
{\bf \S 1. Introduction.} Let $f:z\mapsto z^2+c$ be a quadratic polynomial
with connected filled Julia set $K$. The Julia set $\partial K$ is as usual
denoted by $J$. Let $\psi: \overline{\BBB C}\smallsetminus \overline{\BBB D}
\rightarrow \overline{\BBB C}\smallsetminus K$ be the unique conformal
isomorphism, normalized as $\psi(\infty)=\infty$ and $\psi'(\infty)=1$,
which conjugates the squaring map to $f$: 
\begin{equation} 
\label{eqn:Bottcher} 
\psi(z^2)=f(\psi(z)). 
\end{equation} 
(The inverse $\psi^{-1}$ is often called the {\it B\"{o}ttcher} coordinate.)
By the {\it external ray} $R_t$ we mean the image of the radial line $\{
\psi(re^{2 \pi i t}): r>1 \}$, where $t\in \BBB T=\BBB R /\BBB Z$ is the
{\it angle} of the ray. We say that $R_t$ {\it lands} at $z\in J$ if
$\lim_{r\rightarrow 1} \psi(re^{2 \pi i t})=z$. A point $z\in J$ is called
{\it accessible} if there exists an external ray which lands at $z$, and it
is called {\it biaccessible} if $z$ is the landing point of more than one
ray. It can be shown that $z\in J$ is biaccessible if and only if
$K\smallsetminus \{ z \}$ is disconnected (see \cite{McMullen}, p. 85).  
 
Let us denote by $\gamma(t)$ the radial limit $\lim_{r\rightarrow 1}
\psi(re^{2 \pi i t})$. According to a classical theorem of Fatou (see for
example \cite{Rudin}, p. 249), $\gamma(t)$ exists for almost every $t\in
\BBB T$ in the sense of the Lebesgue measure. For all such angles $t$, it
follows from (\ref{eqn:Bottcher}) that $\gamma$ conjugates the doubling map
to the action of $f$ on the Julia set: 
$$\gamma(2t)=f(\gamma(t)).$$ 
When $K$, or equivalently $J$, is locally-connected, it follows from the
theorem of Carath\'{e}odory that $\gamma$ is defined and continuous on the
whole circle. In this case, the surjective map $\gamma :\BBB T\rightarrow J$
is called the {\it Carath\'{e}odory loop}. Evidently the biaccessible points
in $J$ correspond to the points where $\gamma$ fails to be one-to-one. 
 
Whether or not $J$ is locally-connected, the Lebesgue measure on the circle
$\BBB T$ pushes forward by $\gamma$ to a probability measure $\mu$ on the
Julia set. Complex analysts call $\mu$ the ``harmonic measure'' on $J$, but
in the context of holomorphic dynamics, $\mu$ is called the {\it Brolin
measure}. It has the following nice properties: 
\begin{enumerate} 
\item[(i)] 
The support of $\mu$ is the whole Julia set, with $\mu(J)=1$. 
\item[(ii)] 
$\mu$ is invariant under the $180^{\circ}$ rotation $z\mapsto -z$, i.e.,
$\mu(-A)=\mu(A)$ for every measurable set $A\subset J$. 
\item[(iii)] 
$\mu$ is $f$-invariant, i.e., $\mu(f^{-1}(A))=\mu(A)$ for every measurable
set $A\subset J$. Moreover, $\mu$ is ergodic in the sense that for every
measurable set $A\subset J$ with $f^{-1}(A)=A$, we have $\mu(A)=0$ or
$\mu(A)=1$. 
\end{enumerate} 
All of these properties are immediate consequences of the corresponding
properties of the Lebesgue measure and the angle-doubling map on the unit
circle.  
Properties (ii) and (iii) are equivalent to the next property, which will be
used repeatedly in this paper: 
\begin{enumerate} 
\item[(iv)] 
$\mu(f(A))=2 \mu(A)$ for every measurable set $A\subset J$ for which the
restriction $f|_A$ is one-to-one. 
\end{enumerate} 
Brolin proved that with respect to this measure $\mu$ the backward orbits of
typical points have an asymptotically uniform distribution
\cite{Brolin}. Lyubich has proved that $\mu$ is the unique measure of
maximal entropy $\log 2$. He has also constructed such invariant measures of
maximal entropy for arbitrary rational maps of the Riemann sphere
\cite{Lyubich}. 
  
It follows from general plane topology that the set of points in $J$ which
are the landing points of more than two external rays is at most countable
(see for example \cite{Pommerenke}, p. 36). On the other hand, the number of
rays landing at a point is constant along an orbit, unless the orbit passes
through the critical point. It follows from ergodicity that either
$\mu$-almost every point in the Julia set is the landing point of a unique
ray, or else $\mu$-almost every point is the landing point of exactly two
rays. 
 
As an example, for the {\it Chebyshev polynomial} $z\mapsto z^2-2$, the
Julia set is the closed interval $[-2,2]$ on the real line. Here every point
is the landing point of exactly two rays except for the endpoints $\pm 2$
where unique rays land. There are no other known examples of quadratic Julia
sets with two rays landing at almost every point. In fact, as I heard from
J. Hubbard and later M. Lyubich, it is conjectured that a polynomial Julia
set has this property {\it only if} it is a straight line segment in which
case the map is conjugate to a Chebyshev polynomial, up to sign. In this
paper, we will confirm this conjecture for quadratic Julia sets which are
locally-connected. The second part of this paper \cite{Zakeri}, which is an
expanded version of \cite{Schleicher-Zakeri}, considers the Julia sets of
quadratic polynomials with irrationally indifferent fixed points. By a
completely different method we prove that every biaccessible point in $J$
eventually maps to the critical point in the Siegel case and to the Cremer
fixed point otherwise. As a byproduct, it follows that the set of
biaccessible points in the Julia set has Brolin measure zero. This settles
some cases that are not covered by \thmref{main} of this paper, since in the
Cremer case $J$ is certainly non locally-connected, and in the Siegel case
$J$ may or may not be locally-connected. \vspace{0.1 in}\\ 
{\bf Acknowledgement.} I am indebted to Jack Milnor who suggested the
possibility of such a theorem and generously shared his ideas with me, which
play an important role in this paper.\\ \\     
{\bf \S 2. Basic Definitions.} Let $f:z\mapsto z^2+c$ be a quadratic
polynomial whose filled Julia set $K$ is locally-connected. As usual, the
fixed points of $f$ are denoted by $\alpha$ and $\beta$, where $\beta$ is
the more repelling fixed point. If $\alpha$ is attracting or $\alpha=\beta \
(\Leftrightarrow c=1/4)$, the Julia set of $f$ is a Jordan curve with a
unique external ray landing at every point. Hence there are no biaccessible
points at all and \thmref{main} below is trivially true. {\it So we may as
well assume that $\alpha \neq \beta$ and $\alpha$ is not attracting}. It
follows that either $\alpha \in J$, or else $\alpha$ is the center of a
fixed Siegel disk for $f$. 
 
By an {\it embedded arc} in $K$ we mean any subset of $K$ homeomorphic to
the closed interval $[0,1]\subset \BBB R$. Since $K$ is locally-connected,
for any two points $x,y \in K$ there exists an embedded arc $\eta$ in $K$
which connects $x$ to $y$. If $K$ has no interior so that $J=K$ is full,
then $\eta$ is uniquely determined by the two endpoints $x$ and $y$. If $K$
does have interior, however, there is usually more than a choice for
$\eta$. In what follows, we will show how to choose a canonical embedded arc
between any two points in the filled Julia set. 
 
Suppose that $int(K)$ is non-vacuous. Every component $U$ of this interior
is a bounded Fatou component whose closure $\overline U$ is homeomorphic to
the closed unit disk $\overline {\BBB D}$ since $K$ is
locally-connected. According to Fatou and Sullivan (see for example
\cite{Milnor}), every such component eventually maps to a periodic Fatou
component which is either the immediate basin of attraction of an attracting
periodic point, or an attracting petal for a parabolic periodic point, or a
periodic Siegel disk. We refer to these cases simply as {\it hyperbolic},
{\it parabolic} and {\it Siegel} cases. Note that in the hyperbolic and
parabolic cases the critical point $0$ belongs to a central Fatou component
which we denote by $U_0$. Also by our assumption on the $\alpha$-fixed
point, periodic Fatou components in the hyperbolic and parabolic cases form
a cycle of period $>1$.  
 
Next, we would like to choose a ``center'' $c(U)$ in every bounded Fatou
component $U$ subject only to the following conditions: 
\begin{enumerate} 
\item[\bf{(C1)}] 
$c(-U)=-c(U)$, 
\item[\bf{(C2)}] 
If $U$ contains the critical value $c=f(0)$, then $c(U)=c$, 
\item[\bf{(C3)}]   
If $U$ contains the fixed point $\alpha$, then $c(U)=\alpha$. 
\end{enumerate} 
If follows from {\bf (C1)} that whenever the critical point $0$ belongs to
the Fatou set, then it is the center of the corresponding Fatou component
$U_0$: $c(U_0)=0$. Also {\bf(C3)} corresponds to the case where the
$\alpha$-fixed point is the center of a fixed Siegel disk $U$.  
 
Given any bounded Fatou component $U$, there exists a homeomorphism $\phi:
\overline U \iso \overline {\BBB D}$ which is holomorphic in $U$ with
$\phi(c(U))=0$. An arc in $\overline U$ of the form $\phi^{-1}\{ re^{i
\theta}: 0\leq a\leq r\leq b\leq 1 \}$ is called a {\it radial arc}. Since
$\phi$ is unique up to postcomposition with a rigid rotation of
$\overline{\BBB D}$, radial arcs in $\overline U$ are well-defined.      
 
Following \cite{Douady-Hubbard}, we call an embedded arc $I$ in $K$ {\it
regulated} if for every bounded Fatou component $U$, the intersection $I
\cap \overline U$ is either empty or a point or consists of radial arcs in
$\overline U$ (see also \cite{Douady}, where he uses the word ``legal'' for
regulated). 
 
\begin{lem} 
\label{regulated} 
Given any two points $x,y\in K$, there exists a unique regulated arc $I$ in
$K$ with endpoints $x,y$. Furthermore, if $\eta$ is any embedded arc in $K$
which connects $x$ to $y$, then $I \cap J \subset \eta \cap J$.  
\end{lem} 
 
\begin{pf} 
Take any embedded arc $\eta$ in $K$ with endpoints $x,y$. It is easy to see
how one can deform $\eta$ to a regulated arc $I$. Let $U$ be a bounded Fatou
component whose closure intersects $\eta$. Choose any parametrization
$h:[0,1]\rightarrow K$ with $\eta=h([0,1])$, and define  
$$t_0=\inf \{ t\in [0,1]: h(t)\in \overline U \} ,$$ 
$$t_1=\sup \{ t\in [0,1]: h(t)\in \overline U \} .$$ 
In other words, $t_0$ is the first moment $\eta$ hits $\overline U$ and
$t_1$ is the last moment $\eta$ stays in $\overline U$. If $t_0\neq t_1$,
replace the subarc of $\eta$ from $h(t_0)$ to $h(t_1)$ by the radial arc
from $h(t_0)$ to $c(U)$ followed by the radial arc from $c(U)$ to $h(t_1)$
(see \figref{reg}). If $h(t_0)$ and $h(t_1)$ happen to be on the same radial
arc, simply connect the two by the radial arc between them.  
\realfig{reg}{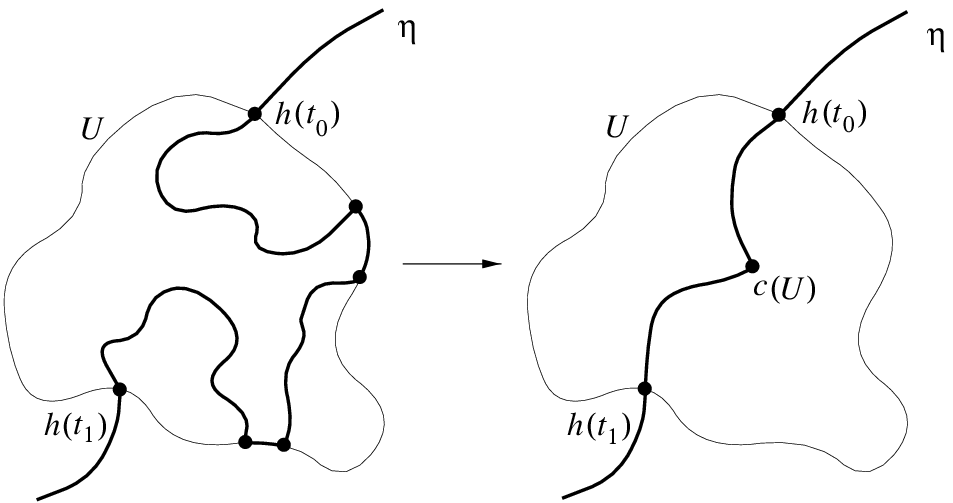}{Deforming an embedded arc to a regulated
arc.}{10cm} 
 
Applying this construction to the intersection with every such Fatou
component, we obtain a regulated arc $I$ with endpoints $x,y$. Evidently we
have the inclusion $I \cap J \subset \eta \cap J$.  
 
To prove uniqueness, suppose that $I$ and $I'$ are both regulated, with the
same endpoints $x,y$. If $I \neq I'$, then the complement $\BBB C
\smallsetminus (I \cup I')$ has a bounded connected component $V$. By the
Maximum Principle, $V$ is contained in some bounded Fatou component $U$. It
follows that the boundary $\partial V$ must be contained in a union of at
most four radial arcs in $\overline U$. But a finite union of radial arcs
cannot bound an open set in $\overline U$. Therefore, $I=I'$. 
\end{pf} 
 
The regulated arc $I$ given by the above lemma is denoted by $[x,y]$. The
open arc $(x,y)$ is defined by $[x,y]\smallsetminus \{ x, y\}$, and
similarly we can define the semi-open arc $[x,y)$. 
 
More generally, given finitely many points $x,y, \ldots, z$ in $K$, there is
a unique smallest connected set $\tree \subset K$ made up of regulated arcs
which contains all of these points. In fact this set is always a (finite)
topological tree. We call $\tree$ the {\it regulated tree} generated by $\{
x,y,\ldots, z\}$. A vertex of this tree with exactly one edge attached to it
is called an {\it end} of the tree. A point which is not an end is called an
{\it interior point} of the tree. It follows easily from {\bf(C1)} that  
\begin{equation} 
\label{eqn:sym} 
[-x, -y, \ldots, -z]=- \tree 
\end{equation} 
 
In the case of three distinct points, $[x,y,z]$ is either homeomorphic to a
closed interval or to a letter Y. The first case occurs if and only if one
of the points belongs to the regulated arc connecting the other two. In the
second case, the three points $x,y,z$ are ends of the tree $[x,y,z]$. In
other words, there is a unique interior point $p\in [x,y,z]$ such that
$[x,p] \cap [y,p] = [x,p] \cap [z,p] = [y,p] \cap [z,p] =\{ p \} $ (see
\figref{1}). In this case, we call $[x,y,z]$ a {\it tripod}. Point $p$ is
called the {\it joint} of this tripod. 
\realfig{1}{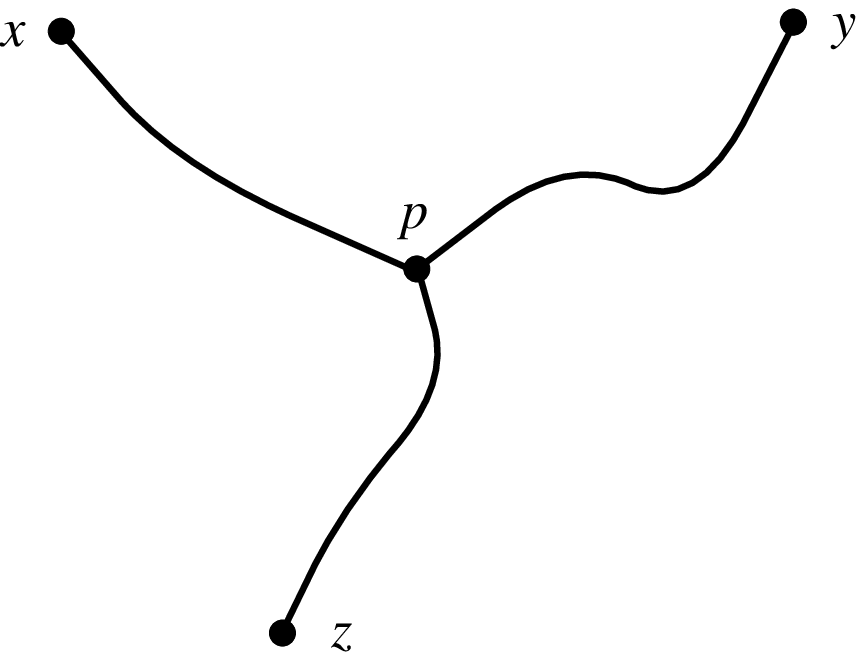}{A tripod $[x,y,z]$ with joint $p$.}{5cm} 
 
The regulated trees as defined above are not preserved by the dynamics of
$f$. In fact, when $K$ has interior, the center of a bounded Fatou component
$U$ is not necessarily mapped by $f$ to that of $f(U)$. Hence regulated arcs
in $U$ do not map to regulated arcs in $f(U)$. This difficulty can be most
conveniently overcome by deforming the polynomial $f$ rel the Julia set into
a new map $F$ which respects the centers. To this end, it suffices to note
that for every bounded Fatou component $U$, there is a homeomorphism between
$\overline{U}$ and the cone over $\partial U$ which sends $c(U)$ to the cone
point and restricts to the identity map on $\partial U$. We can define $F$
so as to preserve this cone structure on various bounded Fatou
components. For example, for any component $U$ and any $p \in \partial U$
take the Poincar\'{e} geodesic in $U$ between $c(U)$ and $p$ and define
$F:\overline{U} \rightarrow f(\overline{U})$ so as to map this geodesic
isometrically to the unique Poincar\'{e} geodesic between $c(f(U))$ and
$f(p)\in \partial f(U)$. (Note that by our assumption $f(U)\neq U$ unless
$U$ is a fixed Siegel disk for which the $\alpha$ fixed point is the
center. So in any case $\alpha$ is still a fixed point of $F$.) Apply this
construction to every bounded Fatou component and let $F=f$ anywhere
else. The map $F$ will be the required modification of $f$ which satisfies
the following properties:  
\begin{enumerate} 
\item[\bf{(F1)}] 
$F(c(U))=c(F(U))$ for every bounded Fatou component $U$. In particular, by
{\bf(C2)}, whether or not the critical point $0$ belongs to the Fatou set,
$F(0)=f(0)=c$ is always the critical value of $f$. 
\item[\bf{(F2)}] 
$F=f$ on the closure of the basin of attraction of infinity. 
\item[\bf{(F3)}] 
$F(z)=F(z') \Leftrightarrow z=\pm z'$. 
\item[\bf{(F4)}] 
$\alpha$ and $\beta$ are the only fixed points of $F$.  
\end{enumerate} 
Also, since the support of the Brolin measure is the Julia set where $f$ and
$F$ agree, it follows that properties (iii)-(iv) in section {\bf \S 1} also
hold for $F$. In other words, 
\begin{enumerate} 
\item[\bf{(F5)}] 
$\mu(F^{-1}(A))=\mu(A)$ for any measurable set $A\subset \BBB C$, and 
\item[\bf{(F6)}] 
$\mu(F(A))=2\mu(A)$ for any measurable set $A\subset \BBB C$ for which $
F|_A$ is one-to-one. 
\end{enumerate} 
 
\begin{lem} 
\label{respect} 
Let $x,y, \ldots, z\in K$. Suppose that the critical point $0$ is not an
interior point of the tree $\tree$. Then $F$ maps $\tree$ homeomorphically
to $[F(x), F(y),\ldots, F(z)]$. 
\end{lem} 
 
In this case, we simply write  
$$F:[x,y,\ldots, z] \iso [F(x), F(y),\ldots, F(z)].$$ 
 
\begin{pf} 
First let us show that $F$ restricted to $\tree$ is injective. If not, it
follows from {\bf(F3)} that $\tree$ contains a pair $\pm a$ of symmetric
points. By (\ref{eqn:sym}), we see that $[a,-a]=-[a,-a]$. Hence the
$180^{\circ}$ rotation from the arc $[a,-a]$ to itself must have a fixed
point, namely the critcal point $0$. But this implies that $0$ is an
interior point of $\tree$, contrary to our assumption.   
 
Therefore, $F$ restricted to $\tree$ is injective. The image tree $F(\tree)$
is evidently connected and contains all of the image points
$F(x),F(y),\ldots, F(z)$. Since all the ends of $F(\tree)$ are among
$F(x),F(y),\ldots, F(z)$, we conclude that it is also minimal. To finish the
proof, it is enough to show that the image of every regulated arc in $\tree$
is a regulated arc. But this follows from {\bf(F1)} since $F$ preserves the
centers hence the radial arcs in bounded Fatou components of $f$. 
\end{pf} 
\noindent 
{\bf Definition.} By the {\it spine} of the filled Julia set $K$ we mean the
unique regulated arc $[-\beta, \beta]$ between the $\beta$-fixed point and
its preimage $-\beta$, which are the landing points of the unique external
rays $R_0$ and $R_{1/2}$ respectively. By (\ref{eqn:sym}), the spine is
invariant under the $180^{\circ}$ rotation $z\mapsto -z$. In particular, the
critical point $0$ always belongs to the spine.\\ 
 
Let $z\in J$ be a biaccessible point, with a ray pair $(R_t, R_s)$ landing
at $z$ and $0<t<s<1$. If $z\notin [-\beta, \beta]$, it follows that both $t$
and $s$ satisfy $0<t<s<1/2$ or $1/2<t<s<1$. Consider the orbit of the ray
pair $(R_t, R_s)$ under $f$. Since there exists an integer $n > 0$ such that
$1/2 \leq  2^n s-2^n t < 1$, the corresponding rays $f^{\circ n}(R_t)$ and
$f^{\circ n}(R_s)$ must belong to different sides of the curve $R_{1/2}\cup
[-\beta, \beta] \cup R_0$ (see \figref{one}). Therefore $f^{\circ n}(z)\in
[-\beta, \beta]$. This means that the set $B$ of all biaccessible points in
the Julia set is contained in the union of preimages of the spine: 
\begin{equation} 
\label{eqn:BB} 
B\subset \bigcup_{n\geq 0}f^{-n}[-\beta, \beta]. 
\end{equation} 
 
\realfig{one}{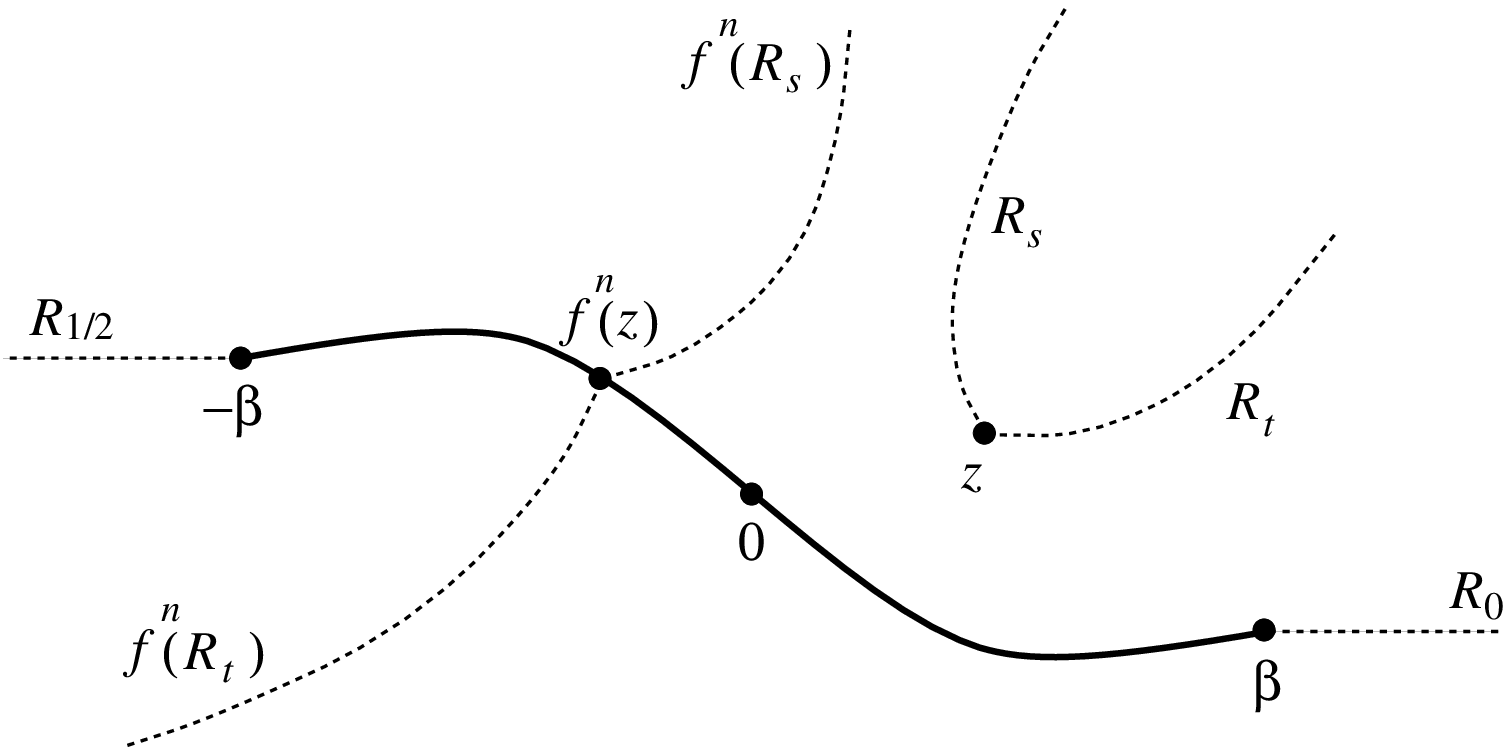}{}{9cm} 
\vspace{0.15in} 
\noindent 
{\bf \S 3. Main Theorem and Supporting Lemmas.} In this paper we will prove
the following theorem: 
 
\begin{thm} 
\label{main} 
If the Julia set $J$ of the quadratic polynomial $f:z\mapsto z^2+c$ is
locally-connected, then the set of all biaccessible points in $J$ has Brolin
measure zero unless $f$ is the Chebyshev polynomial $z\mapsto z^2-2$ for
which the corresponding measure is one.  
\end{thm} 
 
By (\ref{eqn:BB}), it suffices to show that for every non-Chebyshev
quadratic, the Brolin measure $\mu[-\beta, \beta]$ of the spine is zero.  
 
The proof depends on several lemmas which will be given in this section and
section {\bf \S 4}. 
 
\begin{lem} 
\label{arc} 
 
\noindent 
\begin{enumerate} 
\item[(a)] 
Any point in the Julia set $J$ which belongs to the boundary of two Fatou
components is necessarily biaccessible. 
\item[(b)] 
Let $\eta $ be any embedded arc in the filled Julia set $K$ and $z$ be a
point in $\eta \cap J$ which is not an endpoint of $\eta$. Then either $z$
is biaccessible or it belongs to the boundary of a unique bounded Fatou
component.  
\end{enumerate} 
\end{lem} 
 
\begin{pf} 
(a) Let $U$ and $U'$ be two such Fatou components, with $z \in \partial U
\cap \partial U'$. Assume that $z$ is not biaccessible. Then $K
\smallsetminus \{ z \}$ is connected, so there exists an embedded arc $\eta
$ in $K$ between $c(U)$ and $c(U')$ which avoids $z$. By \lemref{regulated},
$I \cap J \subset \eta \cap J$, where $I=[c(U),c(U')]=[c(U),z]\cup [z,
c(U')]$ is the unique regulated arc between $c(U)$ and $c(U')$. It follows
that $\eta$ must contain $z$, which is a contradiction. 
  
(b) If $z$ is not biaccessible, then $K \smallsetminus \{ z \}$ is
connected. Hence there exists an embedded arc $\eta'$ in $K$ between the two
endpoints of $\eta$ which avoids $z$. Take a bounded connected component $V$
of the complement $\BBB C \smallsetminus (\eta \cup \eta')$ which contains
$z$ in its closure. By the Maximum Principle, $V$ must be contained in a
bounded Fatou component. Hence $z$ belongs to the boundary of this bounded
Fatou component. Uniqueness follows from part (a). 
\end{pf} 
 
\begin{cor} 
\label{b-b} 
Let $f:z\mapsto z^2+c$ have locally-connected Julia set. If the
$\alpha$-fixed point is not attracting and $\alpha \neq \beta$, then neither
the $\beta$-fixed point nor any of its preimages can belong to the boundary
of a bounded Fatou component of $f$.  
\end{cor} 
 
\begin{pf} 
Otherwise there exists a bounded Fatou component $U$ with $\beta \in
\partial U$. Hence $\beta \in \partial U \cap \partial f(U)$. If $U=f(U)$,
it must be a fixed Siegel disk by the assumption. But in this case
$f|_{\partial U}$ is conjugate to an irrational rotation so it cannot have a
fixed point. Therefore $U\neq f(U)$. By \lemref{arc}(a), $\beta$ will be
biaccessible. But this is impossible since the $\beta$-fixed point is always
the landing point of the unique ray $R_0$. 
\end{pf} 
\vspace{0.15in} 
\noindent  
{\bf Remark.} In the non locally-connected case, it is {\it not} known if
the $\beta$-fixed point can be on the boundary of any bounded Fatou
component. In fact, it is not known if there are examples of quadratic
polynomials with a fixed Siegel disk whose boundary is the whole Julia
set. Any such quadratic would provide a counterexample to the above
corollary in the non locally-connected case.    
 
\begin{lem} 
\label{bet} 
If $x\notin [-\beta, \beta]$, then $[-\beta, x, \beta]$ is a tripod. 
\end{lem} 
 
\begin{pf} 
Otherwise, we must have $-\beta\in (x, \beta)$ or $\beta\in (x, -\beta)$. In
either case, it follows that $-\beta$ or $\beta$ belongs to the interior of
an embedded arc in the filled Julia set. But $\beta$ is the landing point of
the unique ray $R_0$. Since the orbit $-\beta \mapsto \beta$ does not pass
through the critical point, it follows that $-\beta$ is also the landing
point of the unique ray $R_{1/2}$. By \lemref{arc}(b), either $\beta$ or
$-\beta$ must be on the boundary of a bounded Fatou component, which
contradicts \corref{b-b}.  
\end{pf} 
 
Here is a definition which will be used repeatedly in all subsequent
arguments:\\ \\ 
{\bf Definition.} We define a projection $\pi:K\rightarrow [-\beta, \beta]$
as follows: For $x\in [-\beta, \beta]$, let $\pi(x)=x$. If $x\notin [-\beta,
\beta]$, then $[-\beta, x, \beta]$ is a tripod by \lemref{bet}, and we
define $\pi(x)\in (-\beta, \beta)$ to be the joint of this tripod.\\ 
 
Note that $\pi(x)$ can be described as the unique point in $[-\beta, \beta]$
such that for any $y$ on the spine, $[x,\pi(x)]\subset [x,y]$. Set
theoretically $\pi$ a retraction from $K$ onto its spine. However, when $K$
has interior, $\pi$ is not continuous. 
 
For simplicity, we denote the regulated arc $[x, \pi(x)]$ by $I_x$. Since
$\pi(-x)=-\pi(x)$, we have $I_{-x}=-I_x$. 
 
\begin{lem} 
\label{AL} 
The $\alpha$-fixed point belongs to $(-\beta, 0)$.  
\end{lem} 
 
\begin{pf} 
First we prove that $\alpha\in (-\beta, \beta)$. In fact, if $\alpha$
belonged to $J$ and were off the spine, then the external rays which land at
$\alpha$ would all belong to one side of the curve $R_{1/2}\cup [-\beta,
\beta] \cup R_0$. This would contradict the fact that the angle-doubling map
on the circle has no forward orbit which is entirely contained in the
interval $(0,1/2)$ or $(1/2,1)$. On the other hand, if $\alpha$ belonged to
the Fatou set and were off the spine, then it would have to be the center of
a fixed Siegel disk whose closure by {\bf(C3)} touches  $[-\beta, \beta]$ at
the unique point $0$. Take the external ray $R_t$ which lands at the
critical value $c$. Since the entire orbit of $c$ is on one side of the
curve $R_{1/2}\cup [-\beta, \beta] \cup R_0$, the forward orbit of $t$ under
the doubling map must be entirely contained in one of the intervals
$(0,1/2)$ or $(1/2,1)$, which is again a contradiction. Therefore,
$\alpha\in (-\beta, \beta)$.    
 
Now suppose that $\alpha\in (0,\beta)$. Then $[\alpha, \beta] \subset
(0,\beta]$. Hence $F:[\alpha, \beta]\iso [\alpha, \beta] $ by
\lemref{respect}. By {\bf(F4)}, there is no fixed point of $F$ in $(\alpha,
\beta)$. Suppose that $[\alpha, \beta] \subset J$. Then $f$ repels all
points in $[\alpha, \beta]$ close to $\alpha$ and $\beta$. Since $f=F$ on
the Julia set, the same must be true for $F$. Hence there has to be an
attracting fixed point for $F$ somewhere in $(\alpha, \beta)$, which is a
contradiction. Therefore $[\alpha, \beta]$ intersects a bounded Fatou
component $U$. Passing to some iterate $f^{\circ n}(U)=F^{\circ n}(U)$, we
may as well assume that $U$ is periodic. Since $F$ acts monotonically on
$[\alpha, \beta]$, $U$ must be fixed. Hence $U$ is a Siegel disk with
$c(U)=\alpha$. Now $\partial U$ intersects $[\alpha, \beta]$ at a unique
point $p$ which is not the $\beta$-fixed point by \corref{b-b}. Clearly
$F(p)=p$, which is a contradiction. This shows that  $\alpha \in (-\beta,
0)$, and completes the proof. 
\end{pf} 
 
\begin{lem} 
\label{OM} 
There exists an $F$-preimage $\omega$ of $0$ in $(-\beta, \alpha)$. The
other preimage $-\omega$ is then in $(-\alpha, \beta)$. 
\end{lem} 
 
\begin{pf} 
$F:[-\beta, \alpha]\iso [\beta, \alpha]$ by \lemref{respect} since $0\notin
(-\beta, \alpha)$ by \lemref{AL}. Again by \lemref{AL} we have $0\in (\beta,
\alpha)$, which shows there exists a unique $\omega\in (-\beta, \alpha)$
with $F(\omega)=0$. 
\end{pf} 
 
Let $<$ be the natural order between the points of the spine induced by any
homeomorphism $h:[-\beta, \beta]\iso [-1,1]\subset \BBB R$, with $h(\pm
\beta)=\pm 1$. In other words, for $x,y\in [-\beta, \beta]$ we have $x<y$ if
and only if $h(x)<h(y)$.  
\begin{cor} 
Let $\pm \omega$ be the two $F$-preimages of $0$ as in \lemref{OM}. Then we
have the following order between the points on the spine (see
\figref{omega}): 
$$-\beta < \omega < \alpha < 0 < -\alpha < -\omega <\beta.$$ 
\end{cor}  
\realfig{omega}{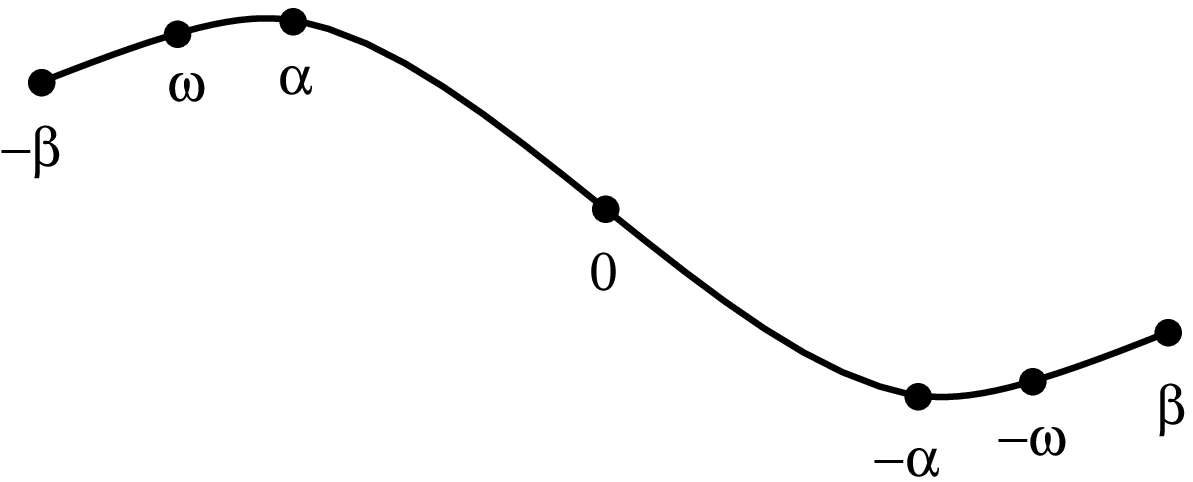}{}{7cm} 
 
\begin{lem} 
\label{CV} 
Let $c=f(0)=F(0)$ be the critical value. Then $\pi(c)\in [-\beta,
\alpha]$. If $\pi(c)=-\beta$, then $c=-\beta$ in which case $f(z)=z^2-2$.  
\end{lem} 
 
\begin{pf} 
By \lemref{respect} we have $F:[0, \beta]\iso [c, \beta]=I_c \cup
[\pi(c),\beta]$. Since $-\alpha \in [0, \beta]$, by {\bf(F3)} and {\bf(F4)}
we must have $F(-\alpha)=\alpha \in [c, \beta]$. This is possible only if
$\alpha \in [\pi(c), \beta]$, which is equivalent to $\pi(c)\in [-\beta,
\alpha]$ (see \figref{2}). 
 
If $\pi(c)=-\beta$, then $c=-\beta$ by \lemref{bet}. It is easy to see that
$z\mapsto z^2-2$ is the only quadratic polynomial with the critical orbit
$0\mapsto c \mapsto \beta$.  
\end{pf} 
\realfig{2}{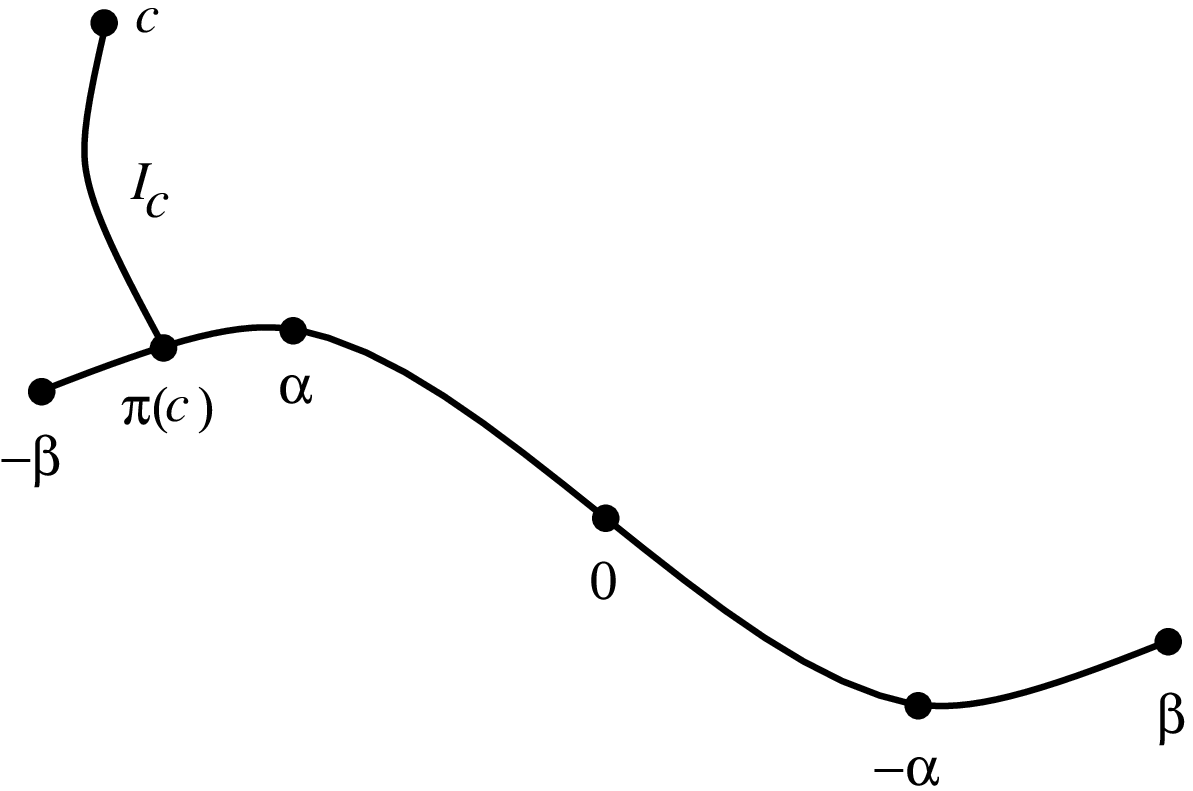}{}{7cm} 
 
\begin{lem} 
\label{XI} 
Suppose that $f$ is not the Chebyshev polynomial. Let
$f(\xi)=F(\xi)=-\beta$. Then $\xi$ does not belong to the spine $[-\beta,
\beta]$. Furthermore, $\pi(\xi) \in [-\alpha, \alpha]$ and
$F(\pi(\xi))=\pi(c)$, with  
$$c\in [-\beta, \beta] \Leftrightarrow \pi(\xi) = 0.$$  
\end{lem} 
 
\begin{pf} 
First suppose that $ \pi(\xi) \neq 0$. Replacing $\xi$ by $-\xi$ if
necessary, we may assume that $\pi(\xi) \in (0, \beta)$. Then $F :
[\xi,\beta] \iso [-\beta, \beta]$, hence $-\alpha \in [\xi, \beta]$ which
implies that $-\alpha\in [\pi(\xi),\beta]$, or equivalently, $\pi(\xi)\in
(0, -\alpha]$. Also, since $0 \notin [\xi, \beta]$, $c$ cannot belong to the
spine $[-\beta, \beta]$. By \lemref{CV}, $\pi(c)\in (-\beta, \alpha]$. By
\lemref{respect} the set $[ \xi,0, \beta]$ maps homeomorphically to the
tripod $[-\beta,c, \beta]$, hence it must also be a tripod, with $\xi \notin
[-\beta, \beta]$, and with the joint $\pi(\xi)$ mapped to $\pi(c)$ by $F$
(see \figref{3}). 
 
Now suppose that $ \pi(\xi)=0$. Then by a similar argument, the set $[\xi,
0,  \beta]=[\xi, \beta]$ still maps homeomorphically to the spine $[-\beta,
\beta]$ since it does not contain a pair of symmetric points about the
origin. In particular, $c$ must belong to the spine. By \lemref{CV},
$c=\pi(c) \in (-\beta, \alpha]$.   
\end{pf} 
\realfig{3}{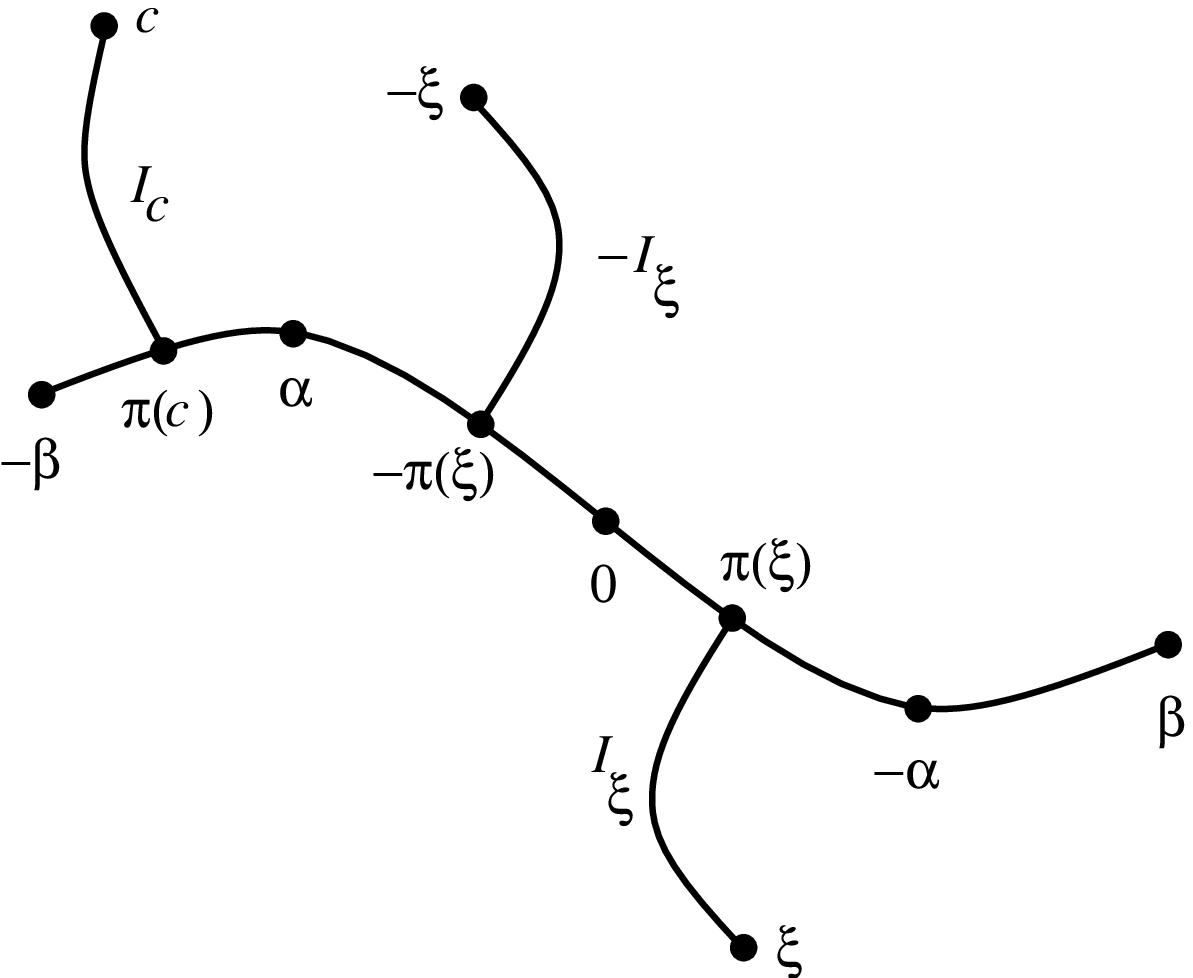}{}{7cm} 
 
\begin{cor} 
\label{map} 
$F$ maps $[0, \pm \pi(\xi)]$ to $I_c$ and $\pm I_{\xi}$ to $[-\beta,
\pi(c)]$ homeomorphically (see \figref{3}). 
\end{cor} 
 
Thus in all non-Chebyshev cases we have the situation illustrated in
\figref{3} (except that $I_c$ may collapse to a point $\Leftrightarrow$
$[-\pi(\xi), \pi(\xi)]$ may collapse to a point, or alternatively $\pi(c)$
may coincide with $\alpha$). Here 
$$\pm \xi \stackrel{F}{\mapsto} -\beta \stackrel{F}{\mapsto} \beta,$$ 
$$\pm \pi(\xi) \stackrel{F}{\mapsto} \pi(c),$$ 
and 
$$\pm \omega \stackrel{F}{\mapsto} 0 \stackrel{F}{\mapsto} c,$$ 
where $\omega$ lies somewhere between $-\beta$ and $\alpha$. 
 
\begin{lem} 
\label{enough} 
Suppose that $f$ is not the Chebyshev polynomial. Then the Brolin measure
$\mu[-\beta, \beta]$ of the spine is zero if and only if $\mu(I_c)=0$.  
\end{lem} 
 
Note that the condition $\mu(I_c)=0$ is trivially satisfied if $c=\pi(c)$
belongs to the spine. The latter happens, for example, when the Julia set of
$f(z)=z^2+c$ with $c\in \BBB R$ is full. When the Julia set is full, it is
conjectured that the critical value belongs to the spine if and only if $c$
is real.  
 
\begin{pf}  
By \lemref{XI}, for one preimage $\xi$ of $-\beta$, we have $\pi(\xi)\in [0,
-\alpha]$, and then the other preimage $-\xi$ satisfies $\pi(-\xi)\in
[\alpha, 0]$. For simplicity, let $z_0=\pi(\xi)$ and $z_n=F^{\circ
n}(z_0)$. It follows from \corref{map} that 
\begin{equation} 
\label{eqn:preimage} 
F^{-1}([-\beta, \beta]\cup I_c)=[-\beta, \beta]\cup I_{\xi} \cup -I_{\xi}. 
\end{equation} 
By property (ii) in section {\bf \S 1} and {\bf(F5)}, we have 
\begin{equation} 
\label{eqn:equal1}  
\mu (I_{\xi})=\mu (-I_{\xi})=\frac{1}{2} \mu(I_c). 
\end{equation} 
Note that $z_1 =\pi(c) \in (-\beta, \alpha]$ by \lemref{XI} and
\lemref{CV}. By \corref{map}, {\bf(F6)} and (\ref{eqn:equal1}), 
\begin{equation} 
\label{eqn:equal2} 
\mu [-\beta, z_1]=\mu ( F(I_{\xi}))=2\mu(I_{\xi})=\mu(I_c). 
\end{equation}  
If $\mu[-\beta, \beta]=0$, then $\mu [-\beta, z_1]=0$, hence $\mu(I_c)=0$ by
(\ref{eqn:equal2}). Conversely, if $\mu(I_c)=0$, then $\mu [-\beta, z_1]=0$.

To prove $\mu [-\beta, \beta]=0$, we distinguish two cases:\\ 
 
$\bullet${\it Case 1.} $z_1 \in [\omega, \alpha]$. Then $\mu [-\beta,
\omega]\leq \mu [-\beta, z_1]=0$. Hence $\mu [0, \beta]=2 \mu [-\beta,
\omega] =0$, which by symmetry implies $\mu [-\beta, \beta]=0$.\\ 
 
$\bullet${\it Case 2.} $z_1 \in (-\beta, \omega)$. Then $z_2=F(z_1)\in
F(-\beta, \omega)=(0, \beta)$ and $\mu [z_2, \beta]=2 \mu [-\beta,
z_1]=0$. If $z_2 \in [0, -\omega]$, then $\mu [-\omega, \beta]=0$ and it
follows by an argument similar to {\it Case 1} that $\mu [-\beta,
\beta]=0$. So let us assume that $z_2 \in (-\omega, \beta)$. We can repeat
the above argument by considering $z_3=F(z_2)\in (0, \beta)$. If $z_3 \in
[0, -\omega]$, we have $\mu [-\beta, \beta]=0$, otherwise $z_3 \in (-\omega,
\beta)$ and we continue. If this process never stops, it follows that $z_n
\in (-\omega, \beta)$ and $(z_{n+1}, \beta] \supset [z_n , \beta]$ for all
$n$. The limit of the monotone sequence $\{ z_n \}$ will then be a fixed
point of $F$ in $(-\omega, \beta)$, which contradicts {\bf(F4)}. 
\end{pf} 
\vspace{0.15 in} 
\noindent 
{\bf \S 4. The Proof.} The idea of the proof of \thmref{main} is as follows:
We consider the $n$-th iterate of $c=f(0)=F(0)$,  $c_n=F^{\circ
n}(c)$. Under the assumption $\mu(I_c)>0$, we show that $c_n$ cannot belong
to the spine and the Brolin measure of the arc $I_{c_n}$ tends to infinity
as $n\rightarrow \infty$, which is clearly impossible since
$\mu(J)=1$. Hence we must have  $\mu(I_c)=0$. By \lemref{enough}, this
proves the theorem.\\ \\ 
{\bf Definition.} Let $I_1$ and $I_2$ be two regulated arcs in the filled
Julia set $K$. We say that $I_1$ and $I_2$ {\it overlap} if the intersection
$I_1\cap I_2$ contains more than one point. It follows that $I_1\cap I_2$ is
a nondegenerate regulated arc $I$ in $K$. We often say that $I_1$ and $I_2$
{\it overlap along} $I$.\\  
 
It is not hard to check that for $x,y \in K\smallsetminus [-\beta, \beta]$,
the arcs $I_x$ and $I_y$ overlap if and only if $x$ and $y$ belong to the
same connected component of $K\smallsetminus [-\beta, \beta]$. In
particular, we must have $\pi(x)=\pi(y)$. 
 
\begin{lem} 
\label{3case} 
Let $x\in K\smallsetminus [-\beta, \beta]$. Then one and only one of the
following cases occurs, as illustrated in Figures 7, 8, 9: 
 
\begin{enumerate} 
\item[(a)] 
$I_x$ and $I_{\xi}$ (or $-I_{\xi}$) overlap along an arc $I_y$. Then $F$
maps $[x,y]$ homeomorphically to $I_{F(x)}=[F(x),F(y)]$. 
\item[(b)] 
$\pi(x) \in (-\pi(\xi), \pi(\xi))$. Then $F$ maps $I_x$ homeomorphically to
the arc $F(I_x)=[F(x), F(\pi(x))]$. In this case, $I_{F(x)}$ and $I_c$
overlap along $I_{F(\pi(x))}=[F(\pi(x)), \pi(c)]$.  
\item[(c)] 
$\pi(x) \notin (-\pi(\xi), \pi(\xi))$ and $I_x$ and $\pm I_{\xi}$ do not
overlap. Then $F$ maps $I_x$ homeomorphically to $I_{F(x)}$. 
\end{enumerate} 
\end{lem} 
 
\begin{pf} 
(a) If $x\in I_{\xi}$ or $-I_{\xi}$, then $y=x$ and the result it
trivial. Otherwise, $[\xi, x, \pi(x)=\pm \pi(\xi)]$ maps homeomorphically to
$[-\beta, F(x), \pi(c)]$ (see \figref{7}). Hence $F(y)=\pi(F(x))$ and the
result follows. 
 
(b) If $\pi(x) \in (-\pi(\xi), \pi(\xi))$, then $F(\pi(x))\in
I_c\smallsetminus \{ \pi(c) \}$, hence $I_{F(x)}$ and $I_c$ overlap along
$I_{F(\pi(x))}$ (see \figref{6}). 
 
(c) Since $\pi(x) \notin (-\pi(\xi), \pi(\xi))$, $F(\pi(x))\in [-\beta,
\beta]$. So the claim is proved once we show that $\pi(F(x))=F(\pi(x))$. If
these two points are distinct, then the nondegenerate arc
$I=[\pi(F(x)),F(\pi(x))]\subset [-\beta, \beta]$ is contained in $[F(x),
F(\pi(x))]$ (see \figref{5}). Hence $F^{-1}(I)$ will be a nondegenerate arc
in $I_x \cap I_{\xi}$ or $I_x \cap -I_{\xi}$, which contradicts our
assumption. 
\end{pf} 
 
\begin{figure}[htp]
  \hbox to \hsize{
   \hbox{\vbox{\hsize=.45\hsize \columnwidth=\hsize
	  \centerline{\psfig{figure=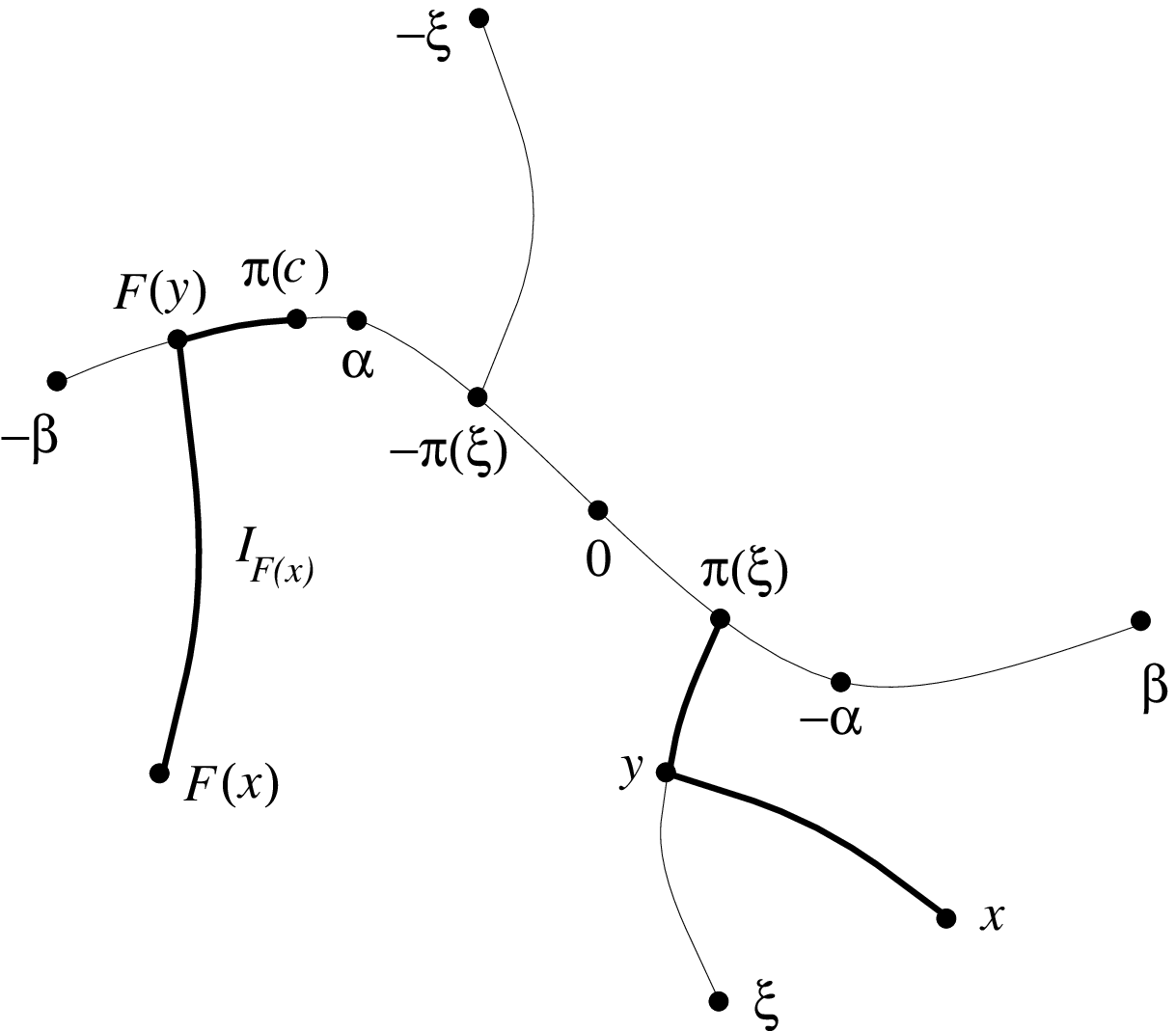,width=6cm}}
          \caption[7]{}
	  \oplabel{7}
 	}}
    \hfil
    \hbox{\vbox{\hsize=.45\hsize \columnwidth=\hsize
	  \centerline{\psfig{figure=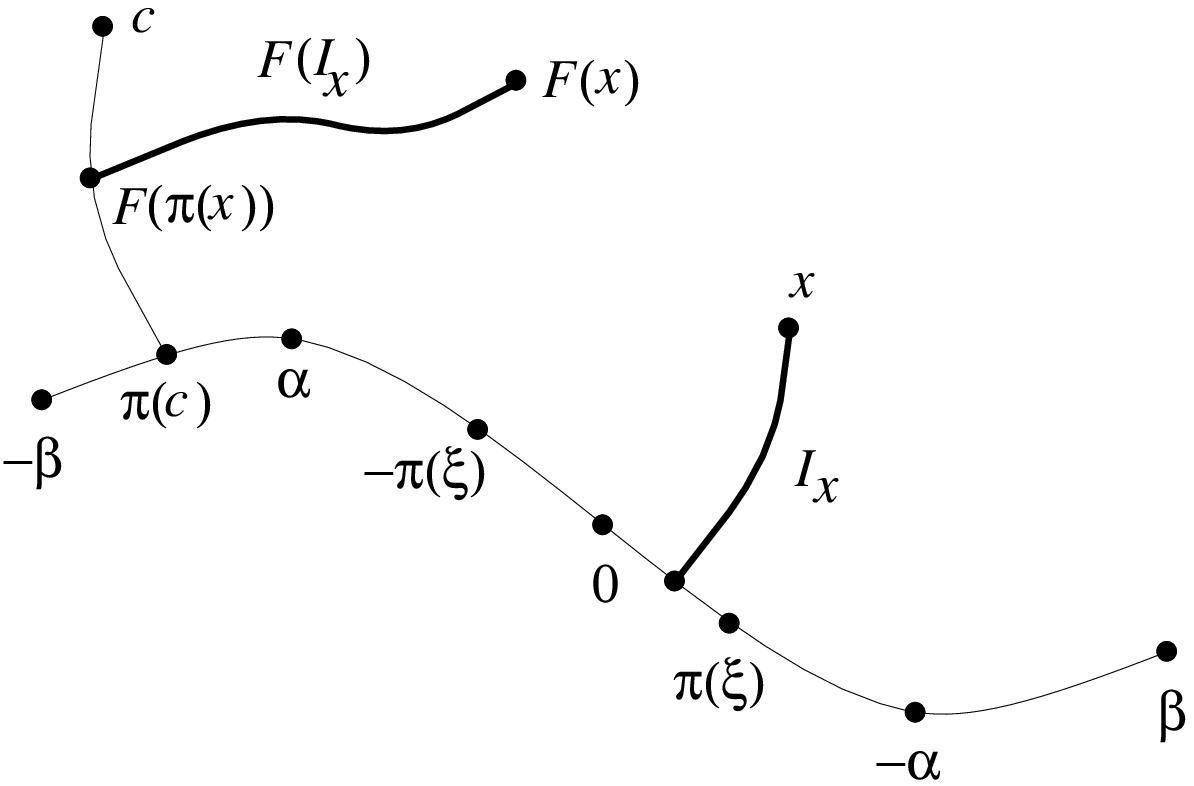,width=6cm}}
	  \caption[6]{}
	  \oplabel{6}
         }}
   }
   \vbox{
	  \centerline{\psfig{figure=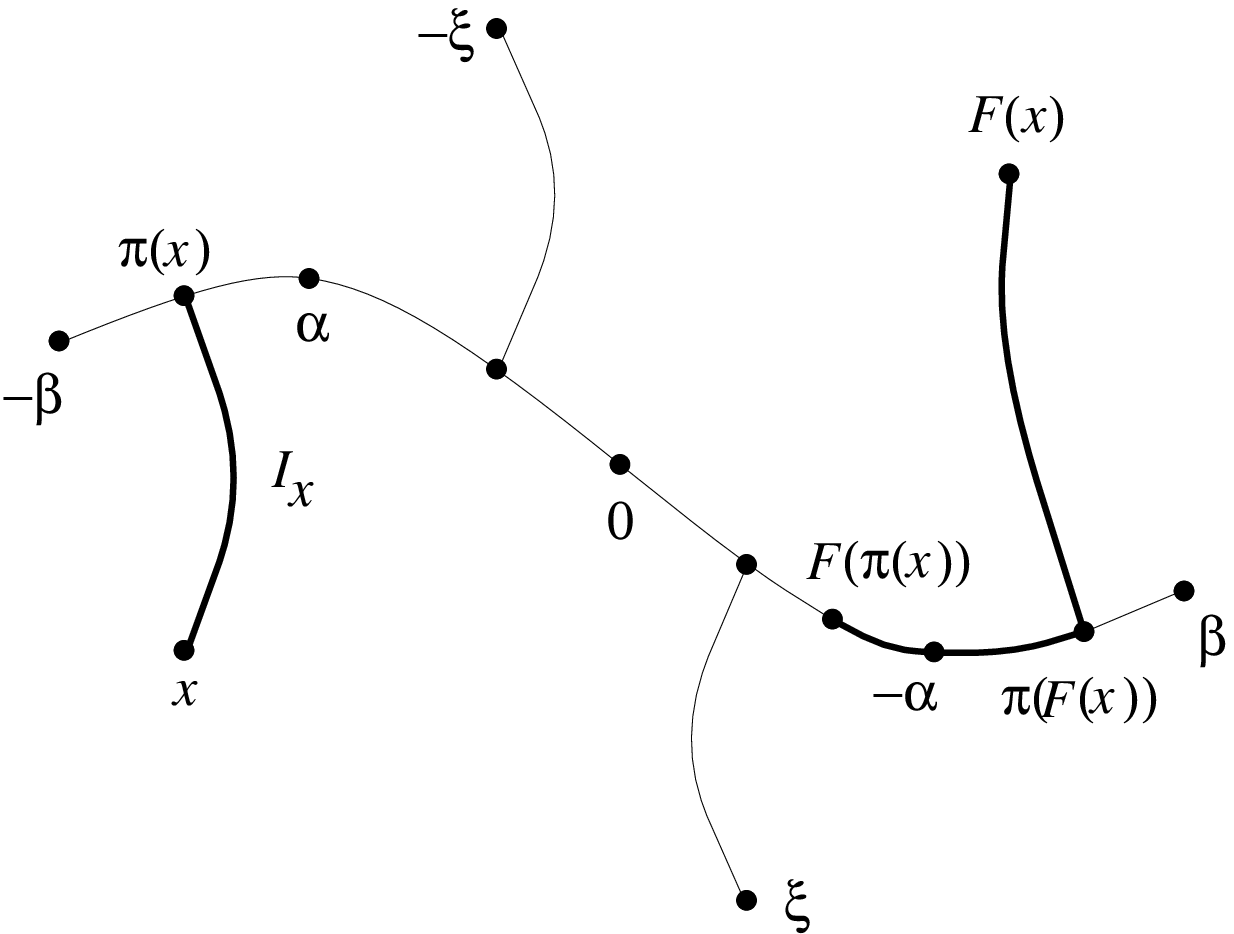,width=6cm}}
	  \caption[5]{}
	  \oplabel{5}
         }
\end{figure}

\goodbreak
Let us put $m=\mu(I_c)$. By (\ref{eqn:equal1}), we have $\mu(\pm
I_{\xi})=m/2$. 
  
\begin{cor} 
\label{oncebigger} 
If $x\in K \smallsetminus [-\beta, \beta]$ and $\mu (I_x)\geq 2m$, then
$\mu(I_{F(x)})\geq \displaystyle{ \frac{3}{2}} \mu(I_x)$. 
\end{cor} 
 
\begin{pf} 
By \lemref{3case} one and only one of the cases (a)-(c) occurs. In case (b),
we have $\mu(I_{F(x)})=\mu(F(I_x))+\mu(I_{F(\pi(x))})\geq
\mu(F(I_x))=2\mu(I_x)$ and in case (c), $\mu(I_{F(x)})=2\mu(I_x)$. In case
(a), 
$$\begin{array}{rl} 
\mu(I_{F(x)}) & = \mu[F(x),F(y)]=2\mu[x,y]\\ 
              & = 2(\mu(I_x)-\mu(I_y))\\ 
              & \geq 2(\mu(I_x)-\mu(\pm I_{\xi}))\\ 
              & = 2\mu(I_x)-m\\ 
              & \geq (3/2)\mu(I_x), 
\end{array}$$ 
which proves the corollary. 
\end{pf} 
 
\noindent 
{\it Proof of \thmref{main}.} Consider the orbit of the critical value $\{
c=c_0, c_1, c_2, \ldots \} $, where $c_n=F^{\circ n}(c)$. Let
$m=\mu(I_c)>0$, and apply \lemref{3case} to the point $x=c$. Clearly the
only possible cases are (a) and (c), since $\pi(c)\notin (-\pi(\xi),
\pi(\xi))$.  
 
In case (c) we obtain the estimate $\mu(I_{c_1})\geq 2m$. This, by repeated
application of \corref{oncebigger}, will lead to the estimate
$\mu(I_{c_{n+1}})\geq (3/2)^n \mu(I_{c_1})$ which tends to infinity as
$n\rightarrow \infty$ and therefore is impossible.  
 
In case (a), $I_c$ and $-I_{\xi}$ overlap along some $I_y$ with $F(y)\in
(-\beta, \pi(c))$ and $\mu(I_{c_1})=\mu[c_1,F(y)] = 2\mu[c,y] \geq m$. Apply
\lemref{3case} this time to $x=c_1$. Note that the only possible case is
(c), since $\pi(c_1)\in [-\beta, \pi(c))$. This gives the estimate
$\mu(I_{c_2})=2\mu(I_{c_1})\geq 2m$. Hence successive applications of
\corref{oncebigger} will give the estimate $\mu(I_{c_{n+1}})\geq (3/2)^n
\mu(I_{c_2})$, which again contradicts the fact that the Brolin measure of
the Julia set is finite. 
 
The contradiction shows that $m=\mu(I_c)$ must be zero, and this completes
the proof of \thmref{main} by \lemref{enough}.\ \ \hfill$\Box$\\  
 
It has been shown recently that the Julia set of a {\it real} quadratic
polynomial in the Mandelbrot set is locally-connected \cite{Levin}. These
correspond to quadratics $f:z\mapsto z^2+c$ with $-2\leq c \leq
1/4$. Therefore, we have the following corollary of \thmref{main}: 
 
\begin{cor} 
Let $f:z\mapsto z^2+c$ be a real quadratic polynomial with $-2< c \leq
1/4$. Then, with respect to the Brolin measure on the Julia set of $f$,
almost every point is the landing point of a unique external ray. 
\end{cor}   
\vspace{0.15 in} 
\noindent 
{\bf \S 5. Further Discussion.} Finally, we consider the following result,
which is a consequence of \thmref{main} as well as the fact that the Julia
set has no compact forward-invariant proper subsets of positive Brolin
measure. 
 
\begin{thm} 
\label{invisible} 
Let $f:z\mapsto z^2+c$ be a quadratic polynomial with locally-connected
filled Julia set $K$. If we exclude the Chebyshev case and the cases where
the $\alpha$-fixed point of $f$ is attracting or $\alpha=\beta$, then every
embedded arc in $K$ has Brolin measure zero. 
\end{thm} 
 
The exceptional cases correspond respectively to $c=-2$ where the Julia set
is a straight line segment, $c$ in the ``main cardioid'' of the Mandelbrot
set where the Julia set is a quasicircle, and $c=1/4$ where the Julia set is
a Jordan curve but not a quasicircle. Roughly speaking, the theorem says
that in any other case, embedded arcs are buried in the filled Julia set so
that they are almost invisible from the basin of infinity. 
 
We need the following elementary observation for the proof: 
     
\begin{lem} 
\label{density} 
Let $A\subset J$ be forward-invariant under $f$, i.e., $f(A)\subset A$. Then
either $\mu(A)=0$ or $\mu(A)=1$. In particular, if $A$ is compact and $A\neq
J$, then $\mu(A)=0$. 
\end{lem} 
\begin{pf} 
Let $\gamma :\BBB T\rightarrow J$ be the Carath\'{e}odory loop and
$E=\gamma^{-1}(A)$. Then $E$ is forward-invariant under the doubling map
$d:\BBB T\rightarrow \BBB T$ defined by $d(t)=2t$ (mod 1). We prove that
$\ell(E)=0$ or $\ell(E)=1$, where $\ell$ denotes the Lebesgue measure on
$\BBB T$. Let $\ell(E)>0$ and let $x$ be a point of density of $E$. Given an
$\varepsilon >0$, we can find an $n>0$ and an interval $S\subset \BBB T$
centered at $x$ such that $\ell(S)=2^{-n}$ and $\ell(S\cap E)\geq
(1-\varepsilon) \ell(S)$. Apply the $n$-th iterate $d^{\circ n}$ on $S$ and
use $d^{\circ n}(E)\subset E$ to estimate 
$$1-\varepsilon \leq 2^n \ell (S\cap E)=\ell (d^{\circ n}(S\cap E))\leq \ell
(\BBB T\cap E)=\ell (E).$$ 
Since this is true for every $\varepsilon >0$, we must have $\ell(E)=1$. 
\end{pf} 
 
\begin{cor} 
\label{boundary} 
Still assuming that $K$ is locally-connected, the Brolin measure of the
union of the boundaries of bounded Fatou components of $f$ is zero unless
the $\alpha$-fixed point is attracting or $\alpha=\beta$ in which case the
corresponding measure is one. 
\end{cor} 
 
\begin{pf} 
Since every bounded Fatou component eventually enters a cycle of Fatou
components of the form $U_1 \mapsto U_2 \mapsto \ldots \mapsto U_p \mapsto
U_1$, it suffices to prove that $\mu(A)=0$, where
$A=\bigcup_{j=1}^{p}\partial U_j$. This set is compact and forward-invariant
under $f$, so by \lemref{density} if $\mu(A)>0$, then $A=J$ must be the
case. But this implies that $f$ has only $p$ bounded Fatou components. It is
easy to see that this can happen only if $p=1$, in which case the component
is either the immediate basin of attraction for an attracting fixed point or
the attracting petal for a parabolic fixed point. 
\end{pf} 
 
As an illustrative example, consider a quadratic polynomial $f$ whose
$\alpha$-fixed point is the center of a Siegel disk $U$ with rotation number
$\theta$ of constant type (an example is provided by $f:z\mapsto
z^2-0.3905408-0.5867879 i$, where $\theta =(\sqrt{5}-1)/2$ is the golden
mean). By \cite{Petersen}, the filled Julia set is locally-connected. The
critical point $0\in \partial U$ is the landing point of exactly two rays
$(R_s, R_{s+1/2})$, where 
$$s=\sum_{0 < p/q < \theta} 2^{-(q+1)}.$$ 
Since the orbit of $0$ is dense on $\partial U$, the set of angles $t$ for
which $\gamma (t) \in \partial U$ coincides with the closure of the orbit of
$s$ under the doubling map on the circle. This set is known to be an
invariant Cantor set $C$ of measure zero in the interval $[s, s+1/2]\subset
\BBB T$ (see \cite{Bullett}). It follows that the set of all $t$ for which
$\gamma (t)$ belongs to the boundary of a bounded Fatou component is the
countable union of Cantor sets consisting of $C$ and all its preimages under
the doubling map. This set has Lebesgue measure zero, hence the union of the
boundaries of all bounded Fatou components will have Brolin measure zero. \\
\\ 
{\it Proof of \thmref{invisible}.} Let $\eta \subset K$ be any emdedded
arc. Let $B$ be the set of biaccessible points in $J$ and $B'$ be the set of
all points in $J$ which belong to the boundary of a bounded Fatou
component. By \thmref{main} and \corref{boundary}, we have
$\mu(B)=\mu(B')=0$. On the other hand, by \lemref{arc}(b), every $z\in \eta
\cap J$ is either an endpoint or it belongs to $B\cup B'$. 
Hence, $\mu(\eta)=\mu(\eta \cap J)\leq \mu(B\cup B')=0$.\ \  \hfill$\Box$ 
 
\begin{cor} 
A locally-connected quadratic Julia set is not a countable union of embedded
arcs unless it is a straight line or a Jordan curve. 
\end{cor}

\end{document}